\documentclass[12pt,amsfonts]{article}
\usepackage{amsmath,amsthm, amssymb, amstext}
\begin{document}

\title{Scalars convected by a 2D incompressible flow}

\author{Diego Cordoba \\
{\small Department of Mathematics} \\
{\small University of Chicago} \\
{\small 5734 University Av, Il 60637} \\
{\small Telephone: 773 702-9787, e-mail: dcg@math.uchicago.edu} \\
{\small and} \\
Charles Fefferman\thanks{Partially supported by NSF grant DMS 0070692.}\\ 
{\small Princeton University} \\
{\small Fine Hall, Washington Road, NJ 08544} \\
{\small Phone: 609-258 4205, e-mail: cf@math.princeton.edu} \\
}

\date{January 17 2001}
\maketitle

\markboth{2D incompressible flows}{D.Cordoba and C.Fefferman}

\newtheorem {Thm}{Theorem}
\newtheorem {Def}{Definition}
\newtheorem {Lm}{Lemma}
\newtheorem {prop}{Proposition}
\newtheorem {Rem}{Remark} 
\newtheorem {Cor}{Corollary}
\def\cal{\mathcal}
\newtheorem {Ack*}{Acknowledgments}

\section{Abstract}
    We provide a test for numerical simulations, for several two dimensional incompressible flows, that appear to develop sharp fronts. We show that in order to have a front the velocity has to have uncontrolled velocity growth.

\section{Introduction}
The aim of this paper is to study the possible formation of sharp fronts in finite time for a scalar convected by a two dimensional divergence-free velocity field, with $x=(x_1,x_2) \in R^2$ or $R^2/Z^2$, and $t\in[0,T)$ with $T\leq \infty$. The  scalar function $\theta(x,t)$ and the velocity field $u(x,t) = (u_1(x,t), u_2(x,t)) \in R^2$ satisfy the following set of equations 
\begin{eqnarray}
\left (\partial_t + u\cdot \nabla \right )\theta = 0 \\
\nabla^{\perp}\psi = u \nonumber, 
\end{eqnarray}
where $\nabla_{x}^{\perp} f = (-\frac{\partial f}{\partial x_2}, \frac{\partial f}{\partial x_1})$ for  scalar functions f. The function $\psi$ is the stream function.

There are many physical examples where the solutions satisfy the equations above, with an extra equation or operator that relates $\theta$ with the velocity field. Examples include; Passive scalars, Unsteady Prandtl equations, 2D incompressible Euler equations, Boussinesq, 2D Ideal Magnetohydrodynamics and the Quasi-geostrophic equation.

  In the literature on numerical simulations for the 2D Ideal Magneto-hydrodynamics (MHD) a standard candidate for a current sheet formation (see Fig. 1) is when the level sets of the magnetic stream function ( represented in (1) by $\theta$)  contain a hyperbolic saddle (an X-point configuration). The front is formed when the hyperbolic saddle closes, and becomes two Y-points configuration joined by a current sheet. (See Parker~\cite{parker:94}, Priest-Titov-Rickard~\cite{priest-titov:95}, Friedel-Grauer-Marliani~\cite{friedel-grauer-etal:1997} and Cordoba-Marliani~\cite{cordoba-marliani:1999}.)

 \setlength{\unitlength}{1mm}
\begin{picture}(140,44)(10,-10)
 \linethickness{1pt}
 \put(70,5){\line(1,1){10}}
 \put(70,25){\line(1,-1){10}}
 \put(110,15){\line(1,1){10}}
 \put(110,15){\line(1,-1){10}}
 \put(80,15){\line(1,0){30}}
 \qbezier(49,-1)(25,15)(1,-1)
 \qbezier(49,31)(25,15)(1,31)
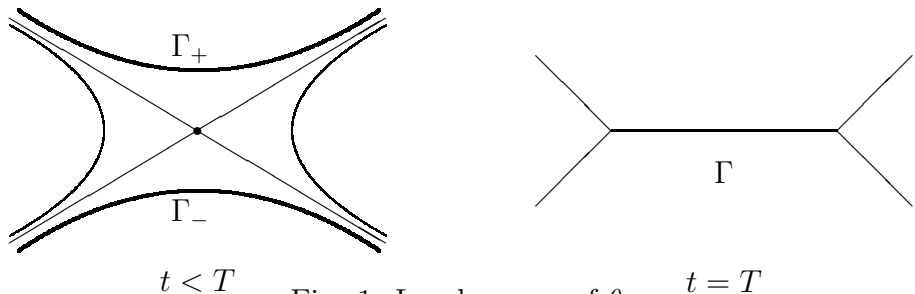
 \put(60,-8){\makebox(0,0){Fig. 1. Level curves of $\theta$.}}
 \put(25,-5){\makebox(0,0){$t < T$}}
 \put(95,-5){\makebox(0,0){$t = T$}}
 \thinlines
 \qbezier(50,1)(25,15)(50,29)
 \qbezier(0,1)(25,15)(0,29)
 \put(0,0){\line(5,3){50}}
 \put(0,30){\line(5,-3){50}}
 \put(25,15){\circle*{1}}
 \put(24,24){\makebox(0,0)[b]{$\Gamma_+$}}
 \put(24,2){\makebox(0,0)[b]{$\Gamma_-$}}
 \put(95,8){\makebox(0,0)[b]{$\Gamma$}}
\end{picture}

  The same configuration was observed in numerical simulations for the Quasi-geostrophic equation (QG). In this case the geometry of the level sets of the temperature has a hyperbolic structure (See Constantin-Majda-Tabak~\cite{constantin-majda-etal:1994b}, Okhitani-Yamada~\cite{Ohkitani-Yamada:97}, Cordoba~\cite{cordoba2:1998} and Constantin-Nie-Schorghofer~\cite{constantin-nie-etal:98}). The QG literature discusses X-points, but not Y-points. In the case of Boussinesq there is no mention, on any numerical simulation study, that a possible singularity is due to the closing of a hyperbolic saddle. In the work of Pumir-Siggia~\cite{pumir-siggia:92} there has been observed evidence for a formation of a front in finite time, across which $\theta$ varies dramatically, on a cap of a symmetric rising bubble.  E-Shu~\cite{e-shu:94} performed numerical simulations with the same initial data as in ~\cite{pumir-siggia:92}, which  suggest that the thickness of the bubble decreases only exponentially.  

The equations for MHD, QG and Boussinesq are as follows

$\textbf{MHD}$:
\begin{eqnarray*}
(\partial_t + u \cdot \nabla ) \theta & = & 0 \\
 (\partial_t + u \cdot \nabla ) \omega & = &\nabla^\perp \theta \cdot \nabla (\Delta \theta) \\
u & = & \nabla^\perp \psi
\end{eqnarray*}
and initial conditions $\theta(x,0) = \theta_{0}$ and $u(x,0) = u_0$. The $\nabla^\perp \theta$ represents the magnetic field, $\Delta \theta$ represents the current density and $\omega = - \Delta\psi $ the vorticity.

$\textbf{QG}$:
\begin{eqnarray*}
(\partial_t + u \cdot \nabla ) \theta = 0 
\end{eqnarray*}
$$
u = \nabla^{\perp}\psi \quad where \quad
\theta=-(-\triangle)^{\frac{1}{2}}\psi
$$
and initial condition $\theta(x,0) = \theta_{0}$. The temperature is represented by $\theta$.

$\textbf{Boussinesq}$:
\begin{eqnarray*}  
\left (\partial_t + u\cdot\nabla_x \right ) \theta & = & 0 \\ 
\left (\partial_t + u\cdot\nabla_x \right ) \omega & = & -\theta_{x_1} \\
u & = & \nabla^\perp \psi 
\end{eqnarray*} 
Again, $\theta$ and u are specified at time t=0. 

\section{Criterion}

A singularity can be formed by collision of two particle trajectories. A trajectory X(q,t) is obtain by solving the following ordinary differential equation 
\begin{eqnarray*} 
\frac{d X(q,t)}{dt} & = & u(X(q,t),t)\\
X(q,0) & = &  q
\end{eqnarray*}
Therefore,
$$
(X(q,t) - X(p,t))_t \leq |X(q,t) - X(p,t)| |\nabla u |_{L^{\infty}}
$$
$$
|X(q,t) - X(p,t)| \geq|X(q,0) - X(p,0)| e^{-\int_{0}^{t}|\nabla u |_{L^{\infty}}ds}
$$
By this trivial argument; in order to have a collision  the quantity $\int_{0}^{t}|\nabla u |_{L^{\infty}}ds$ has to diverge.

A classic criterion for formation of singularities in fluid flows is the theorem of Beale-Kato-Majda (BKM); (see ~\cite{beale-kato-etal:1984}), which improves the estimate described above, and deals with arbitrary singularities, not just collisions. Analogues of the BKM theorem for the above 2-dimensional equations include the following results

 For $\textbf{MHD}$,a singularity cannot develop at a finite time T, unless we have  
\begin{eqnarray*} 
\int_{0}^{T} {sup_{x}|\omega(x,t)| + sup_{x}|\triangle_x\theta(x,t)|} dt = \infty, 
\end{eqnarray*}  
where $\omega$ denotes the vorticity.(See Caflisch-Klapper-Steele~\cite{Caflisch:1997}.) 
  
 For $\textbf{QG}$, a singularity cannot develop at a finite time T, unless we have  
\begin{eqnarray*} 
\int_{0}^{T} sup_{x}|\nabla_x\theta(x,t)| dt = \infty, 
\end{eqnarray*} 
(See Constantin-Majda-Tabak~\cite{constantin-majda-etal:1994b}).

For $\textbf{Boussinesq}$, if a singularity  develops at a finite time T then  
\begin{eqnarray*} 
\int_{0}^{T} sup_{x}|\omega(x,t)| dt = \infty\ \ and\ \ \int_{0}^{T}\int_{0}^{t} sup_{x}|\nabla_x\theta(x,s)| ds dt = \infty.  
\end{eqnarray*} 
(See E-Shu~\cite{e-shu:94}.) 
 
See also Constantin-Majda-Tabak~\cite{constantin-majda-etal:1994b} and Constantin-Fefferman-Majda~\cite{constantin-fefferman-etal:1996} for other conditions involving direction fields, that rule out formation of singularities in fluids. 

In the case of 2D Euler, a singularity cannot develop at a finite time. From the BKM viewpoint this follows from the fact that $\omega$ is advected by the fluid, and therefore $sup_{x}|\omega(x,t)|$ is independent of t. (See BKM ~\cite{beale-kato-etal:1984}.)
 
 Instead of looking at particle trajectories we look at level curves. Because the scalar function $\theta$ is convected by the flow, that  implies that the level curves are transported by the flow. A possible singular scenario is due to level curves approaching each other very fast which will lead to a fast growth of the gradient of the scalar function. In this paper we present a variant of the BKM criterion for sharp front formation. We provide a test for numerical simulations that appear to develop sharp fronts. The BKM Theorem shows that the vorticity grows large if any singularity forms; our Theorem 1 shows that the velocity grows large if a sharp front forms.

The theorem we present in this paper was announced in~\cite{cordoba-fefferman:01}.

\section{Sharp Fronts}

The scalar function $\theta$ is convected by the flow, therefore the level curves move with the flow. A sharp front forms when two of these level curves collapse on a single curve. We define  two level curves to be two distinct time-dependent arcs $\Gamma_+(t)$, $\Gamma_-(t)$ that move with the fluid and collapse at finite time into a single arc $\Gamma$. More precisely, suppose the arcs are given by

\begin{eqnarray}
\Gamma_{\pm} = \{(x_1,x_2)\in R^2: x_2= f_{\pm}(x,t), x_1\in [a,b]\}\ \ for\ \ 0\leq t < T,\label{eq:1}
\end{eqnarray} 
 with 
\begin{eqnarray}
f_{\pm}\in C^{1}([a,b]\times [0,T)) \label{eq:2}
\end{eqnarray}
and
\begin{eqnarray} 
f_-(x_1,t) < f_+(x_1,t)\ \ for\ \ all\ \ x_1\in[a,b],\ \ t\in[0,T).  \label{eq:3}
\end{eqnarray} 
We call the length b-a of the interval [a,b] the length of the front. The assumption that $\Gamma_{\pm}(t)$ move with the fluid means that  
\begin{eqnarray} 
 u_2(x_1,x_2,t) = \frac{\partial f_{\pm}}{\partial x_1}(x_1,t)\cdot u_1(x_1,x_2,t) + \frac{\partial f_{\pm}}{\partial t}(x_1,t)\ \ at\ \ x_2 = f_{\pm}(x_1,t).\label{eq:4}  
\end{eqnarray} 
This holds in particular for level curves of scalar functions g(x,t) that satisfy $\left (\partial_t + u\cdot\nabla_x \right ) g = 0 $.
The collapse of $\Gamma_{\pm}(t)$ into a single curve $\Gamma$ at time T means here simply that 
\begin{eqnarray} 
 lim_{t\rightarrow T^-}(f_+(x_1,t) - f_-(x_1,t)) = 0\ \ for\ \ all\ \ x_1\in[a,b]. \label{eq:5}
\end{eqnarray} 
and $f_+(x_1,t) - f_-(x_1,t)$ is bounded for all $x_1\in[a,b]$, $ t\in[0,T)$.
 
When (\ref{eq:1}), (\ref{eq:2}), (\ref{eq:3}), (\ref{eq:4}) and (\ref{eq:5}) hold, then we say that the fluid forms a $\textbf{sharp front}$ at time T. 
 
The standard candidates for a singularity for MHD and QG are described by the definition given for a sharp front.  We investigate the possible formation of a sharp front. 
 
 The following  assumption will allow us to rule out formation of sharp fronts.  
  We say that the fluid has $\textbf{controlled velocity growth}$ if we have 
\begin{eqnarray} 
\int_{0}^{T} sup\{ |u(x_1,x_2,t)|: x_1\in[a,b], f_-(x_1,t) \leq x_2 \leq f_+(x_1,t)\}  dt <\infty. \label{eq:6}
\end{eqnarray} 
If (\ref{eq:6}) fails, then we say that the fluid has $\textbf{uncontrolled velocity growth}$.

\begin{Lm} Let $\theta$ be a smooth solution of Eq.1 defined for $t\in[0,T)$. Assume there is a $\textbf{sharp front}$ at time T. Then
\end{Lm}
\begin{eqnarray}
\left (\frac{d}{dt}\right )\left (\int_a^b [f_+(x_1,t) - f_-(x_1,t)] dx_1 \right ) & = & \psi(a,f_+(a,t),t) - \psi(a,f_-(a,t),t)\nonumber \\
& + & \psi(b,f_-(b,t),t) - \psi(b,f_+(b,t),t) \label{eq:7}.
\end{eqnarray}
Proof: Take the derivative of the stream function with respect to $x_1$ along an arc $\Gamma_{\pm}(t)$ 
\begin{eqnarray}
\frac{\partial\psi(x_1,f_{\pm}(x_1,t),t)}{\partial x_1} = u_2(x_1,f_{\pm}(x_1,t),t) - \frac{\partial f_{\pm}}{\partial x_1}u_1(x_1,f_{\pm}(x_1,t),t) \label{eq:8}
\end{eqnarray}
by combining (\ref{eq:8}) and (\ref{eq:4}) we obtain

\begin{eqnarray}
\frac{\partial\psi(x_1,f_{\pm}(x_1,t),t)}{\partial x_1} = \frac{\partial f_{\pm}}{\partial t}(x_1,t)\label{eq:9}
\end{eqnarray}
Expression (\ref{eq:7}) follows from integrating (\ref{eq:9}) with respect to $x_1$ between a and b.

\begin{Thm}Let u(x,t) be a divergence-free velocity field, with controlled velocity growth. Then a $\textbf{sharp front}$ cannot develop at time T. 
\end{Thm}
Proof: Assume there is a sharp front at time T. We define
\begin{eqnarray*}
A(t) = \int_{\tilde{a}(t)}^{\tilde{b}(t)} [f_+(x_1,t) - f_-(x_1,t)] dx_1
\end{eqnarray*}
where 
\begin{eqnarray*}
\tilde{a}(t) = a + \int_{t}^{T} sup\{ |u(x_1,x_2,s)|: x_1\in[a,b], f_-(x_1,s) \leq x_2 \leq f_+(x_1,s)\}  ds
\end{eqnarray*}
and
\begin{eqnarray*}
\tilde{b}(t) = b - \int_{t}^{T} sup\{ |u(x_1,x_2,s)|: x_1\in[a,b], f_-(x_1,s) \leq x_2 \leq f_+(x_1,s)\}  ds
\end{eqnarray*}

There is controlled velocity growth, therefore there exists $t^*\in [0,T)$ such that $\tilde{a}(t)\in [a,b]$ and $\tilde{b}(t)\in [a,b]$ for all $t\in [t^*, T)$.

We take the derivative of A(t) with respect to time
\begin{eqnarray*}
\frac{d A(t)}{d t} = sup |u|\cdot \delta(\tilde{b},t) + sup |u|\cdot \delta(\tilde{a},t) +\int_{\tilde{a}(t)}^{\tilde{b}(t)} \frac{\partial}{\partial t}[f_+(x_1,t) - f_-(x_1,t)] dx_1.
\end{eqnarray*}
where  $sup |u| = sup\{ |u(x_1,x_2,t)|: x_1\in[a,b], f_-(x_1,t) \leq x_2 \leq f_+(x_1,t)\}$ and $\delta(z,t)=f_+(z,t) - f_-(z,t)$.

Using the definition of the stream function, the mean value theorem and (\ref{eq:7}), it is easy to check that $ \frac{d A(t)}{d t} > 0$ for $t > t^*$. This contradicts (\ref{eq:5}) by the dominated convergence theorem.

\begin{Ack*} 
 This work was initially supported by the American Institute of Mathematics. 
\end{Ack*}

\end{document}